\documentclass[10pt,reqno]{amsart}

\newcommand{\Em}[1]{\textbf{#1}}
\newcommand{\SigmaP}{{\sf Sigma}}
{\bfseries}{\rmfamily}

\def\Q{{\Bbb Q}}

\begin{document}

\author[R.~Pemantle]{Robin Pemantle}
\address[R.~Pemantle]{
        Department of Mathematics \\
        University of Pennsylvania\\
        209 South 33rd Street\\
        Philadelphia, PA 19096, USA}
\email{pemantle@math.upenn.edu}
\thanks{Research supported in part by NSF grant \# DMS-0103635 and
by a Research Experiences for Undergraduates supplement from NSF}

\author[C.~Schneider]{Carsten Schneider}
\address[C.~Schneider]{
        Research Institute for Symbolic Computation\\
        Johannes Kepler University Linz\\
        Altenberger Str. 69\\
        A--4040 Linz, Austria}
\email{Carsten.Schneider@risc.uni-linz.ac.at}
\thanks{The second author was supported by the SFB-grant F1305 of the Austrian
FWF and by grant P16613-N12 of the Austrian FWF}

\title{When is $0.999\ldots$ equal to~1?}

\maketitle

\begin{abstract}
A doubly infinite sum, numerically evaluated at between
$0.999$ and $1.001$, turns out to have a nice value.
\end{abstract}

\section{Introduction}

The three dots in the title do not refer to an infinite sequence of 9's,
but to digits that are increasingly hard to compute.  The question
is philosophical: how many 9's do we need to see before we start
to believe that the constant we are computing is probably equal to 1?
In our case, the constant was given by the infinite sum:
\begin{equation} \label{eq:sum}
S := \sum_{j,k =1}^\infty \frac{H_j (H_{k+1}-1)}{jk(k+1)(j+k)}
\end{equation}
where $H_j := \sum_{i=1}^j (1/i)$ are the harmonic numbers.

This question in the title is beyond the scope of a mathematics journal.
There are, however, mathematical papers proving identities that were
discovered because numerical computation pointed to a simple answer.
The anecdotal evidence then accumulates in a misleading manner:
when the conjectured identity is false we are less likely ever to know.
One purpose of the present note is to document a case when we were
able to evaluate the constant and it turned out not to equal the
simple guess (in this case,~1).  In fact we will prove:

{\bf Theorem 1:}
\begin{eqnarray}
S & := & \sum_{j,k =1}^\infty \frac{H_j (H_{k+1}-1)}{jk(k+1)(j+k)}
   \nonumber \\[1ex]
& = & - 4 \zeta (2) - 2 \zeta (3) + 4 \zeta (2) \zeta (3) + 2 \zeta (5)
   \label{TheIdentity} \\[1ex]
& = & 0.999222 \ldots \nonumber
\end{eqnarray}
{\it where $\zeta (n) = \sum_{i=1}^\infty i^{-n}$ denotes the
Riemann zeta function.} \hfill \\[2ex]

The second purpose of this note is to demonstrate a piece of
software that, with a little human intervention, can find
(and prove) this sort of identity.

\section{Background}

The sum~(\ref{eq:sum}) arose in a paper giving bounds on the run time of
the simplex algorithm on a polytope known as a Klee-Minty cube~\cite{BP}.
The Klee-Minty cube is an example of great theoretical importance
to analyses of simplex algorithm run times.  The derivation of the
expression~(\ref{eq:sum}) is, however, of no importance here because
$S$ appears in~\cite{BP} only as an upper bound, demonstrably not sharp,
for the leading coefficient, $c$, of the expected run time.  The numerology
of $S$ is therefore unrelated to the physical origins of $c$.

On the other hand, the form of the summation~(\ref{eq:sum}) does
lead one at least to hope that an exact value might be derived. In
principle, any hypergeometric identity, for example, that holds
for a finite indefinite or definite sum, may be automatically
proved via the Wilf-Zeilberger method~\cite{AequalsB}, and in fact
implementations of WZ-type software often can handle summand
expressions of the complexity of~(\ref{eq:sum}). Furthermore, in
many cases the WZ-machinery will not only prove but also find such
an identity, if it exists, given only the left hand side. This, in
general, does not extend to summations over summands involving no
extra parameter: if $S$ were to equal~1, that fact would not
necessarily be automatically detectable. The relatively simple
form of the summand, however, gave us hope that summation tricks
more specific to harmonic series might unlock the problem (the
identity  $\sum 1/(n(n+1)) = 1$ stands as a beacon of hope).
Indeed, several telescoping and re-summation tricks were initially
tried, removing all but one infinite summation in various ways.
These identities were useful in improving our numerical bounds on
$S$.

The numerical bounds we had on $S$ were not all that good.
It should be noted that the harmonic numbers are themselves sums,
so the expression~(\ref{eq:sum}) is really a quadruple sum.  This
makes it perhaps less surprising that our best rigorous bounds were
no closer that $10^{-3}$.  To make a long story short, summing in
one variable, then using exact values for thousands of row and
column sums and an integral approximation for the remaining terms,
led to our best rigorous bounds, namely
$$0.999197 \leq S \leq 1.00093 \, .$$
At this point, although the exact value of $S$ was of no use to us,
we felt embarrassed to publish numerical bounds on a constant that
we suspected was equal to~1.  The authors of~\cite{BP} then consulted
the experts in harmonic summation, who consulted their computers, and
came up with Theorem~1.  The remainder of this note proves Theorem~1.

\section{Solving the Problem}

If a sum such as~(\ref{eq:sum}) has a nice value, it does not
necessarily follow that a truncated version of the sum has a nice
formula. But, if it has a nice closed form, automatic methods
might succeed in finding it by first computing a recurrence and
afterwards searching for solutions of the recurrence in which the
closed form can be represented. Thus we consider a truncated
version of the sum, namely,
$$S(a,b)=\sum_{k=1}^b\frac{H_{k+1}-1}{k(k+1)}
   \sum_{j=1}^a\frac{H_j}{j(j+k)},$$
that is, the upper limits, instead of being infinite, are taken to
be integer variables, $a$ and $b$.  In this case, our optimism is
rewarded: we will be able to simplify the inner sum so that we can
sum a second time if we make some alterations that disappear when
the upper limit goes to infinity.

\subsection*{The inner sum}
In a first step we compute a closed form evaluation of the inner sum
$$h(a,k):=\sum_{j=1}^af(k,j)$$
with $f(k,j):=\frac{H_j}{j(j+k)}$.  Here we follow the summation
principles given in~\cite{Schneider:04b} that are inspired
by~\cite{AequalsB}.  Note that all the computations are carried out
with the summation package \SigmaP\ in the computer algebra system
Mathematica.  The role of the computer here is to produce equations
which we may then rapidly and rigorously verify.

Our first step is to compute for the definite sum $h(a,k)$ the
recurrence relation
\begin{multline}\label{Equ:Recurrence}
k^2\,h(a,k)
-(k+1)(2k+1)\,h(a,k+1)\\
+(k+1)(k+2)\,h(a,k+2)
=\frac{a(a+k+2)-(a+1)(k+1)H_a}{(k+1)(a+k+1)(a+k+2)}
\end{multline}
with a variation of Zeilberger's creative telescoping trick.

The way we find this is to guess that there is some $r$ and
some constants $c_0 , \ldots , c_r$ depending on $k$ but not $j$,
and some function $g(k,j)$ which we assume to have a relatively
simple form, such that a relation
$$\sum_{s=0}^r c_s f(k+s,j) = g(k,j+1) - g(k,j)$$
holds for all $k , j \geq 1$.  We ask \SigmaP\ to find such a
relation for various classes of functions $g$ and for $r = 1, 2,
\ldots$, until we achieve success with $r=2$, $c_0(k)=k^2$,
$c_1(k)=-(k+1)(2k+1)$, $c_2(k)=(k+1)(k+2)$, and $g(k,j)=-\frac{j
H_j+k+j}{(k+j)(k+j+1)}$.  We may then verify the relation
\begin{equation}\label{Equ:CreaEqu}
c_0(k)f(k,j)+c_1(k)f(k+1,j)+c_2(k)f(k+2,j)=g(k,j+1)-g(k,j)
\end{equation}
by polynomial arithmetic and by using the relation
$H_{j+1}=H_j+\frac{1}{j+1}$.  Summing~\eqref{Equ:CreaEqu} over
$k$ in the interval $\{ 1 , \ldots , a \}$ proves~\eqref{Equ:Recurrence}.

Next we are in the position to discover and prove that
\begin{equation}\label{FiniteRepId}
h(a,k)=
\frac{kH_k^2-2H_k+kH^{(2)}_k+2kH^{(2)}_a}{2k^2}-\frac{(kH_a-1)}{k^2}\sum_{i=1}^k\frac{1}{a+i}-\frac{1}{k}\sum_{i=1}^k\frac{1}{i}\sum_{j=1}^i\frac{1}{a+j}
\end{equation}
holds for all $a,k\geq1$; here $H^{(r)}_k=\sum_{i=1}^k\frac{1}{i^r}$
denotes the \Em{generalized harmonic numbers}.

We do this by asking \SigmaP\ to find solutions
to~(\ref{Equ:Recurrence}) for all $k\geq1$ and then to plug in the
initial conditions, which in this case, are to match the values of
$h(a,k)$ for $k = 1, 2$; see~\cite{Schneider:04b} for further
details.  Once \SigmaP\ has found the expression on the right hand
side of~(\ref{FiniteRepId}), it is again a finite exercise in
polynomial arithmetic to verify that it
satisfies~(\ref{Equ:Recurrence}), and that it satisfies the two
initial conditions.  Since $r$ initial conditions uniquely
determine the solution, we have proved~(\ref{FiniteRepId}).

\subsection*{Some terms vanish in the limit}
At some point, if we do not winnow out some terms that will
disappear in the limit, our computer will begin to balk. Luckily
it is easy to identify some terms that will contribute $o(1)$ to
the definite sum as $a$ and $b$ go to infinity, and may therefore
be ignored in the evaluation of $S$.  We have, for example, the
elementary estimates
\begin{align}\label{LimitToZero}
\lim_{a\rightarrow\infty}\frac{1}{k^2}\sum_{i=1}^k\frac{1}{a+i}=0,&&
\lim_{a\rightarrow\infty}\frac{H_a}{k}\sum_{i=1}^k\frac{1}{a+i}=0&&\text{and}&&
\lim_{a\rightarrow\infty}\frac{1}{k}\sum_{i=1}^k\frac{1}{i}\sum_{j=1}^i\frac{1}{a+j}=0.
\end{align}

Hence, if we define
\begin{equation} \label{eq:S'}
S'(a,b):=\sum_{k=1}^b\frac{H_{k+1}-1}{k(k+1)}
   \frac{kH_k^2-2H_k+kH^{(2)}_k+2kH^{(2)}_a}{2k^2},
\end{equation}
we have
$$\lim_{a,b\rightarrow\infty}S'(a,b)=S$$
by~\eqref{FiniteRepId} and~\eqref{LimitToZero}. Summarizing,
problem~\eqref{TheIdentity} from above reduces to find and prove
the identity
\begin{equation}\label{SubProblem}
\lim_{a,b\rightarrow\infty}S'(a,b)
   =-4\zeta(2)-2\zeta(3)+4\zeta(2)\zeta(3)+2\zeta(5).
\end{equation}

\subsection*{The outer sum}
Our next step is guided by the fact that we know a few infinite
sums in which the $n^{th}$ summand is of the form $n^{-c}$ times a
monomial in the harmonic and generalized harmonic numbers $H_n$
and $H_n^{(p)}$.  While we cannot solve the general problem of
summing all such univariate series, it makes sense to attempt to
manipulate things into this form. Thus we ask \SigmaP\ to try to
write the summand of the right hand side of~(\ref{eq:S'}) in the
form $g(a,k) - g(a,k-1)$ where the only infinite sums that appear
in $g$ are of the form described above.  The program obligingly
produces the function $g(a,k) = A(a,k) + B(a,k) + C(a,k)$ where
$A, B$ and $C$ are given by
\begin{eqnarray}
A(a,b) & := & \frac{1}{2{(b+1)}^2}\Big(6H_b+ 4bH_b+4H_b^2
+3bH_b^2+
   H_b^3+bH_b^3-6bH^{(2)}_a \label{eq:A} \\
&& +2H_b H^{(2)}_a+2bH_b H^{(2)}_a-2H^{(2)}_b- 7bH^{(2)}_b+H_bH^{(2)}_b+
   bH_bH^{(2)}_b\Big), \nonumber \\[2ex]
B(a,b) & := & -\frac{2b^2}{{(b+1)}^2}\Big(H^{(2)}_a+H^{(2)}_b\Big)
   \label{eq:B} \\[2ex]
\mbox{ and } && \nonumber \\[2ex]
C(a,b) & := & (H^{(2)}_a-1)\sum_{i=1}^{b}\frac{H_i}{i^2}-\sum_{i=1}^{b}
   \frac{H_i^2}{i^3}+ \frac{1}{2}\sum_{i=1}^{b}\frac{H_i^3}{i^2}
   +\frac{1}{2}\sum_{i=1}^{b}\frac{H_iH^{(2)}_i}{i^2}. \label{eq:C}
\end{eqnarray}
This time the correctness of the result supplied by \SigmaP\ is
verified by polynomial arithmetic and the definition of $H_k$ without
need to appeal to any uniqueness results for solutions to recurrences.

Summing the relation
$$ g(k)-g(k-1)=\frac{H_{k+1}-1}{k(k+1)}\frac{kH_k^2-2H_k+kH^{(2)}_k
   +2kH^{(2)}_a}{2k^2}$$
over $k$ in the interval $\{ 1 , \ldots , b \}$ proves that
\begin{equation}\label{FiniteRepId2}
S'(a,b)=A(a,b)+B(a,b)+C(a,b)
\end{equation}
where $A, B$ and $C$ are as in~(\ref{eq:A})~--~(\ref{eq:C}).

\subsection*{Euler sums and multiple $\boldsymbol{\zeta}$-values}

Now we must make good on our supposition that the form of $g(a,k)$
leads to sums we can evaluate in terms of the zeta function.  The
limits of $A(a,b)$ and $B(a,b)$ are obvious:
\begin{align}\label{Ids:AB}
\lim_{a,b\rightarrow\infty}A(a,b)=0&&\text{and}&&
   \lim_{a,b\rightarrow\infty}B(a,b)=-4\zeta(2).
\end{align}
The sums in $C(a,b)$ are handled in two steps, the first being to
reduce them to multiple $\zeta$-values, which can be done for any
Euler sum.

The \Em{Euler sum} of index $(p_1 \leq \cdots \leq p_k ; q)$, named
after~\cite{euler} (see also~\cite{berndt-euler}), is defined by
\begin{equation} \label{eq:euler}
S_{p_1 , \ldots , p_k ; q} := \sum_{n=1}^\infty
   \frac{H_n^{(p_1)} H_n^{(p_2)} \cdots H_n^{(p_k)}}{n^q} \, .
\end{equation}
Repeated indices $p_i = p_{i+1}$ are allowed, so the
summand may have powers of generalized harmonic numbers
in the numerator.
Each of the summands in $C := \lim_{a,b \to \infty} C(a,b)$ is
an Euler sum.  Define the \Em{multiple $\zeta$-values} by
$$\zeta (a_1 , \ldots , a_k) := \sum_{n_1 > \cdots > n_k}
   \frac{1}{n_1^{a_1} \cdots n_k^{a_k}} \, .$$
The number $k$ is called the \Em{multiplicity} of
the multiple $\zeta$-value and $a_1 + \cdots + a_k$ is known
as the \Em{weight}.  The number $a_1 , \ldots , a_k$ need not
be increasing, but it is required that $a_1 \geq 2$, because
this is the condition for the sum to be finite.

One may expand each $H_n^{(p_j)}$ of~\eqref{eq:euler}, and in
this way the general Euler sum becomes the sum of
$j_1^{-p_1} \cdots j_k^{-p_k} n^{-q}$ over all
$(k+1)$-tuples satisfying $n \geq j_1 , \ldots , j_k$.  Decomposing
according to the set of distinct values among $j_1 , \ldots , j_k , n$
produces an expression for $S_{p_1 , \ldots , p_k ; q}$ that is
a linear combination of multiple $\zeta$-values with weight
$q + p_1 + \cdots + p_k$ and multiplicity at most $k+1$.
In the same way, one sees that any product of Euler sums or
multiple $\zeta$-values is itself a linear combination of multiple
$\zeta$-values.

To illustrate, these ideas, we write $C$ as a sum of products
of Euler sums and then change these, term by term, into multiple
$\zeta$-values.  From the definitions,
$$C = S_{\emptyset;2} S_{1;2} - S_{1;2} - S_{1,1;3}
   + \frac{1}{2} S_{1,1,1;2} + \frac{1}{2} S_{1,2;2} \, .$$
Of course $S_{\emptyset;2} = \zeta (2)$.  Next,
$$S_{1;2} = \sum_{n \geq j \geq 1} n^{-2} j^{-1}
   = \sum_{n \geq 1} n^{-3} + \sum_{n > j \geq 1} n^{-2} j^{-1}
   = \zeta(3) + \zeta(2,1) \, .$$
The remaining values are
\begin{eqnarray*}
S_{1,1;3} & = & 2 \zeta(3,1,1) + \zeta(3,2) + 2 \zeta(4,1) + \zeta(5) \\
S_{1,1,1;2} & = & 6 \zeta (2,1,1,1) + 3 \zeta(2,2,1) + 3 \zeta(2,1,2)\\
   && + 6 \zeta (3,1,1) + 3 \zeta(3,2) + 3 \zeta (4,1) + \zeta (5) \\
S_{1,2;2} & = & \zeta(2,2,1) + \zeta(2,1,2) + \zeta(3,2) + \zeta(4,1)
   + \zeta(5) \, .
\end{eqnarray*}
These are derived in the same way as $S_{1;2}$, for example,
$$S_{1,1;3} = \sum_{n \geq j,k \geq 1} n^{-3} j^{-1} k^{-1}$$
where the triples decompose into $n=j=k, n=j \geq k, n = k \geq j,
n > j = k, n > j > k$ and $n > k > j$.  Putting this all together gives
\begin{eqnarray}
C & = & (\zeta(2) - 1)(\zeta(3) + \zeta(2,1)) + \zeta(3,1,1) + \zeta(3,2)
   + 3 \zeta(2,1,1,1) \label{eq:C2} \\
&& + 2 \zeta(2,2,1) + 2 \zeta(2,1,2) + \zeta(2,3) \, . \nonumber
\end{eqnarray}
The decomposition of products of $\zeta$-values and
multiple $\zeta$-values into linear combinations of multiple
$\zeta$-values is entirely analogous, for instance,
\begin{equation} \label{eq:prod}
\zeta(2) \zeta(3) = \zeta(2,3) + \zeta (3,2) + \zeta (5) \, .
\end{equation}

\subsection*{Reducing to ordinary $\boldsymbol{\zeta}$-values}

Theorem~1 indicates that reduction to single $\zeta$-values
is possible.  This is worthwhile, because the single $\zeta$-values
may be more easily and accurately obtained than may multiple
$\zeta$-values.  For example, Maple computes 500 digits of
$\zeta(5)$ in under 10 seconds.

As it happens, some Euler sums and multiple $\zeta$-values
may be represented as sums of products of single $\zeta$-values.
Others, it appears, may not, though this is not proved and no
algorithm is known for determining which Euler sums or which
multiple $\zeta$-values may be represented in this way.
For large $k$, the number of $\Q$-linearly independent
multiple $\zeta$-values of weight $k$ is conjectured by
Zagier to grow like $1.32^k$~\cite[page~17]{Flajolet:98},
which outstrips the dimension $e^{O(\sqrt{k})}$ of the monomials
$\zeta(k_1) \cdots \zeta(k_m)$ in single $\zeta$-values
of total degree $k= k_1 + \cdots + k_m$.  However, $k =
3,4,5,6,7$ and~9, the dimensions are known to be equal,
and hence every multiple $\zeta$-value of these weights
is a polynomial over $\Q$ in single $\zeta$-values.

Among the known relations, two classes suffice for our problem,
namely the sum and duality relations.  The sum relation of degree
$k$ says that the sum of all multiple $\zeta$-values of fixed
multiplicity $m$ and weight $k$ is equal to $\zeta (k)$. For
$k=3,5$ and $m=2$, these are
\begin{eqnarray*}
\zeta (3) & = & \zeta(2,1) \\
\zeta (5) & = & \zeta (4,1) + \zeta (3,2) + \zeta (2,3) \, .
\end{eqnarray*}
The sum relations were conjectured by Moen~\cite[Section~3]{hoffman},
known in special cases to Euler~\cite{euler}, and proved
by~\cite{granville97}.

Secondly, one has the duality theorem, conjectured
by~\cite{hoffman} and proved by \cite{zagier}.  These imply
equality between $\zeta (p_1 , \ldots , p_k)$ and $\zeta (q_1 ,
\ldots , q_l)$ of equal weights when the index sequences are
related by an operation akin to transposing a partition.  The
exact definition is given in, e.g.,~\cite[page~29]{Flajolet:98}
but all we will need to here is to represent the multiple
$\zeta$-values of weight~5 and multiplicity~3 and~4 in terms of
those of multiplicity~2:
\begin{eqnarray*}
\zeta (2,2,1) & = & \zeta (3,2) \hspace{1.7in} \\
\zeta (2,1,2) & = & \zeta (2,3) \hspace{1.7in} \\
\zeta (3,1,1) & = & \zeta (4,1) \hspace{1.7in} \\
\zeta (2,1,1,1) & = & \zeta (5) \, .
\end{eqnarray*}

Using the duality relations to eliminate all terms of multiplicity
greater than~2 from~(\ref{eq:C2}), and then using~(\ref{eq:prod}),
gives
$$C = (\zeta(2) - 1) (\zeta(3) + \zeta(2,1)) + \zeta(4,1)
   + 3 \zeta(2) \zeta (3) \, .$$
The sum relations now yield
$$C = 4 \zeta(2) \zeta(3) - 2 \zeta (3) + 2 \zeta(5) \, ,$$
which finishes the proof of Theorem~1.

We remark that one more relation is needed to reduce all
multiple $\zeta$-values of weight~5 to polynomials in
single $\zeta$-values, namely
\begin{equation} \label{eq:star}
\zeta(3,2) = 3 \zeta(2) \zeta(3) - \frac{11}{2} \zeta(5) \; ,
\end{equation}
but that our expression for $C$ can be simplified in terms of
single $\zeta$-values without~(\ref{eq:star}). This
identity~(\ref{eq:star}) can be produced by a formula
in~\cite{Borwein:95} which can transform any double $\zeta$-value
with odd weight to a polynomial in single $\zeta$-values.

\section{Acknowledgements}

The authors are grateful to Doron Zeilberger for routing the
problem from its originator to \SigmaP\ and \SigmaP's human operator.
The reduction of the Euler sum to the form in the conclusion
of Theorem~1 was accomplished in a very {\em ad hoc} manner in
a previous draft, and we are indebted to an anonymous referee
for leading us to the Euler sum literature and elucidating the
more methodical approach taken in the present manuscript.

\end{document}